\documentclass[11pt,a4paper]{article}

\usepackage[T1]{fontenc}
\usepackage[utf8]{inputenc}
\IfFileExists{lmodern.sty}{\usepackage{lmodern}}{\usepackage{ae,aecompl}}

\usepackage{amsmath,amssymb,amsthm}
\usepackage{mathtools}

\usepackage[margin=2.7cm]{geometry}
\usepackage[expansion=false]{microtype}
\usepackage{booktabs}
\usepackage{enumitem}

\usepackage{xcolor}

\usepackage{listings}
\definecolor{leankeyword}{RGB}{0,64,128}
\definecolor{leancomment}{RGB}{96,96,96}

\lstdefinelanguage{lean4}{
  keywords={theorem,lemma,def,structure,inductive,where,by,decide,
            import,open,namespace,end,variable,fun,match,with,
            deriving,instance,if,then,else,let,have,show},
  sensitive=true,
  comment=[l]{--},
  morecomment=[s]{/-}{-/},
  string=[b]",
}

\usepackage[colorlinks=true,linkcolor=blue!50!black,citecolor=blue!50!black,urlcolor=blue!50!black]{hyperref}
\usepackage{cite}
\usepackage{needspace}
\usepackage[capitalise,noabbrev]{cleveref}
\hypersetup{pdftitle={A Machine-Checked Proof that Gardam's A2-tilde Lattice Does Not Have Unique Products},
  pdfauthor={Ibrahim Mian and Shayaan Siddique}}

\newcommand{\Z}{\mathbb{Z}}
\newcommand{\Ftwo}{\mathbb{F}_2}
\newcommand{\Gam}{\Gamma}
\newcommand{\Gab}{\Gamma^{\mathrm{ab}}}
\newcommand{\Atwo}{\tilde{A}_2}
\newcommand{\Afour}{A_4}
\newcommand{\Zmod}[1]{\Z/#1}
\newcommand{\Free}{F(a,b)}
\newcommand{\code}[1]{\texttt{#1}}
\newcommand{\axiomClosure}{\{\code{propext},\ \code{Classical.choice},\ \code{Quot.sound}\}}

\newtheorem{theorem}{Theorem}[section]
\theoremstyle{definition}
\newtheorem{definition}[theorem]{Definition}
\newtheorem{example}[theorem]{Example}

\title{A Machine-Checked Proof that Gardam's
       $\tilde{A}_2$ Lattice\\ Does Not Have Unique Products}

\author{Ibrahim Mian \qquad Shayaan Siddique\\[4pt]
  \normalsize Millennium Research\\
  \normalsize\texttt{\{ibby,shayaan\}@millenniumresearch.ai}\\
  \normalsize\texttt{ibrahimnmian@gmail.com}, \texttt{shayaansiddique02@gmail.com}}

\date{September 2026}

\begin{document}
\maketitle

\begin{abstract}
Kaplansky's zero-divisor conjecture asserts that the group ring of a
torsion-free group over a field has no zero divisors. It holds for every group
with the unique-product property, so a counterexample can only come from a
torsion-free group without unique products. In lectures in 2021, Gardam
announced that the torsion-free $\tilde{A}_2$ lattice
$\Gamma = \langle a,b \mid a b a^2 b^{-1} a^2 b^{-2},\ a b^3 a b^4 a^{-1} b\rangle$
does not have unique products and presented it as a new candidate group: it
has property~(T), and the known methods for proving the conjecture do not
apply to it. To our knowledge, no proof of the announcement has been
published. We give a proof checked by the Lean~4 kernel and stated against
Mathlib's \code{UniqueProds} class. The witness is an explicit pair of finite
subsets with $|A| = 32$ and $|B| = 28$ in which each of the $896$ products
coincides with another product. For $658$ products the certificate is an
identity in the free group; the other $238$ certificates are explicit products
of conjugated relators, $970$ conjugates in all, checked by free reduction. A
homomorphism onto $\mathbb{Z}/42$ shows that each pair $(u,v)$ differs from its
partner $(u',v')$ as a pair of group elements, which is all the theorem
requires.
Together with a homomorphism onto the alternating group $A_4$ it also shows
that the listed words are pairwise distinct, so the sets have exactly $32$ and
$28$ elements. The witness and certificates come from an untrusted search
program and are re-checked by Lean. The development uses only the axioms
\code{propext}, \code{Classical.choice} and \code{Quot.sound}, with no
\code{sorry} and no \code{native\_decide}. The mathematical statement is
Gardam's. To our knowledge this is the first verification in a proof assistant
of a unique-product failure in a torsion-free group; torsion-freeness of
$\Gamma$ is taken from Gardam and is not formalized here.
\end{abstract}

\section{Introduction}\label{sec:intro}

Kaplansky's zero-divisor conjecture asserts that for every field $K$ and every
torsion-free group $G$, the group ring $K[G]$ has no zero
divisors~\cite{kaplansky1970}. It remains open; Gardam's 2021
counterexample~\cite{gardam2021}, a nontrivial unit in the group ring
$\Ftwo[P]$ of the Promislow group $P$, concerns the companion \emph{unit}
conjecture.

A group $G$ has the \emph{unique-product property} if for all finite nonempty
$A, B \subseteq G$ some element of $AB$ can be written as $ab$ with $a \in A$,
$b \in B$ in exactly one way. Unique products imply that $K[G]$ has no zero
divisors for every field~$K$~\cite{passman1977}, so the conjecture is only in
question for torsion-free groups \emph{without} unique products. Such groups
exist~\cite{ripssegev1987,promislow1988} but few constructions are known. In his lectures on
the Kaplansky conjectures, Gardam sorted the known examples into two families:
small cancellation
constructions~\cite{ripssegev1987,steenbock2015,gms2015,arzhantsevasteenbock2023}
and groups with small presentations, due to
Promislow~\cite{promislow1988}, Carter~\cite{carter2014} and
Soelberg~\cite{soelberg2018} (see also~\cite{nielsensoelberg2024}). He noted
that the examples of the second family are known to satisfy the zero-divisor
conjecture, so that an explicit counterexample would need a new group without
unique products~\cite{gardamslides1}. In the same lectures he announced such a
group.

\begin{theorem}[Gardam, announced 2021 {\cite{gardamslides1,gardamslides2}}]\label{thm:gardam}
The torsion-free group
\begin{equation}\label{eq:presentation}
  \Gam \;=\; \bigl\langle\, a, b \;\bigm|\;
    a b a^2 b^{-1} a^2 b^{-2},\;\;
    a b^3 a b^4 a^{-1} b \,\bigr\rangle
\end{equation}
does not have the unique-product property.
\end{theorem}

Gardam described $\Gam$ as an arithmetic $\Atwo$ lattice and gave a classifying
space for it built from seven triangles. He recorded that it has Kazhdan's
property~(T) and is not linear over $\mathbb{C}$, so that the known techniques
for proving the zero-divisor conjecture fail, while noting that $\Gam$ does
satisfy the Baum--Connes and Farrell--Jones conjectures. He also remarked that
SAT solving can be used to approach the unique-product
property~\cite{gardamslides2}. We have found no published proof of
\cref{thm:gardam} and no published witness for it.

\paragraph{Scope.}
The statement of \cref{thm:gardam} is due to Gardam. This paper supplies a
proof of the unique-product failure that is checked by the kernel of the
Lean~4 proof assistant~\cite{moura2021lean}, together with the explicit witness
and certificates on which the proof rests. The witness was found by our own
search; we have not seen Gardam's. Torsion-freeness of $\Gam$, its description
as a lattice, and property~(T) are not formalized here and are taken from the
sources cited; the formal theorem is a statement about the group presented by
\eqref{eq:presentation} and does not depend on them. We prove nothing about
zero divisors. \Cref{sec:limits} lists further limitations, and
\cref{sec:ai} states how AI tools were used in this work.

\paragraph{Contributions.}
We prove \code{gamma\_not\_uniqueProds}: the group $\Gam$ of
\eqref{eq:presentation}, constructed as a Mathlib \code{PresentedGroup}, does
not satisfy Mathlib's \code{UniqueProds} class~\cite{mathlib2020}. The
statement uses the library's definition rather than one of our own, and the
development uses only the three standard axioms of Lean; a script checks this
for every theorem in the compiled library (\cref{sec:lean}). We give the
witness explicitly, with a certificate for each of the $896$ product
coincidences and for each disequality that is used
(\cref{sec:certs}, \cref{app:witness}). The witness and certificates were
produced by an untrusted search program and are re-checked by a checker whose
soundness is proved in Lean (\cref{sec:engine}). A group with torsion fails to
have unique products for trivial reasons, so a priority claim is only
meaningful for torsion-free groups: to our knowledge, no failure of the
unique-product property in a torsion-free group had previously been verified
in a proof assistant.

\section{Background}\label{sec:background}

\subsection{Unique products}

\begin{definition}\label{def:up}
A group $G$ has \emph{unique products} if for all finite nonempty
$A, B \subseteq G$ there exists $x \in AB$ with exactly one representation
$x = ab$, $a \in A$, $b \in B$.
\end{definition}

If $\alpha = \sum_a \lambda_a a$ and $\beta = \sum_b \mu_b b$ are nonzero
elements of $K[G]$ with supports $A$ and $B$, then a uniquely represented
product $ab \in AB$ receives the coefficient $\lambda_a \mu_b \neq 0$ in
$\alpha\beta$, and nothing cancels it. Hence unique products imply the
zero-divisor conjecture for $G$ over every field~\cite{passman1977}. Mathlib
contains this implication: \code{UniqueProds} is the hypothesis under which
the library derives the absence of zero divisors in monoid algebras.

The first torsion-free group without unique products is due to Rips and
Segev~\cite{ripssegev1987}. Promislow~\cite{promislow1988} then showed that the
group $P$ contains a $14$-element set $A$ such that $AA$ has no uniquely
represented element. Further examples were given by Carter~\cite{carter2014},
Soelberg~\cite{soelberg2018}, Steenbock~\cite{steenbock2015}, Gruber, Martin
and Steenbock~\cite{gms2015}, and Arzhantseva and
Steenbock~\cite{arzhantsevasteenbock2023}. Nielsen and
Soelberg~\cite{nielsensoelberg2024} study the smallest sets without unique
products in torsion-free groups.

\subsection{The group \texorpdfstring{$\Gamma$}{Gamma}}\label{sec:gamma}

The description of $\Gam$ in this subsection follows Gardam's
slides~\cite{gardamslides1,gardamslides2} and is not formalized in our
development. The group $\Gam$ is a torsion-free arithmetic $\Atwo$ lattice.
Lattices of this type have property~(T)~\cite{cms1994,bhv2008}; for background
on groups acting on $\Atwo$ buildings see~\cite{cmsz1993}. Gardam gave
property~(T), together with the fact that $\Gam$ is not linear over
$\mathbb{C}$, as the reason that the known techniques for the zero-divisor
conjecture fail for $\Gam$.

One instance of this failure can be made explicit. Fisher and
S\'anchez-Peralta~\cite{fishersanchezperalta} embed the group algebra of every
torsion-free virtually compact special group in a division ring, which proves
the zero-divisor conjecture for a class of groups that is large by the work of
Agol and Wise~\cite{agol2013,wise2021}. An infinite group with property~(T) is
not in this class. Indeed, a finite-index subgroup that is the fundamental
group of a compact special cube complex acts freely on the universal cover,
which is a finite-dimensional CAT(0) cube complex; the subgroup inherits
property~(T), so by Niblo and Reeves~\cite{nibloreeves1997} it fixes a point,
and a free action with a fixed point forces the subgroup to be trivial, so
the group would be finite.

We use one elementary computation. Abelianizing \eqref{eq:presentation} sends
the two relators to the exponent vectors $(5,-2)$ and $(1,8)$ in $\Z^2$, and
\begin{equation}\label{eq:ab}
  \Gab \;\cong\; \Z^2 \big/ \langle (5,-2), (1,8) \rangle
       \;\cong\; \Zmod{42},
  \qquad
  \det\begin{pmatrix} 5 & -2\\ 1 & 8\end{pmatrix} = 42 .
\end{equation}
Explicitly, $a \mapsto -8$, $b \mapsto 1$ defines a homomorphism
$\psi \colon \Gam \to \Zmod{42}$: the relators map to $5\cdot(-8) - 2 = -42$
and to $-8 + 8 = 0$. It is onto because $\psi(b) = 1$.

\subsection{What a certificate of failure must contain}\label{sec:whatcert}

Let $A$ and $B$ be the sets of elements of $\Gam$ represented by two finite
lists of words. To show that no element of $AB$ is uniquely represented, it
suffices to give, for each pair of listed words $(u,v)$, a second pair
$(u',v')$ of listed words such that
\begin{enumerate}[leftmargin=2em]
\item $uv = u'v'$ in $\Gam$, and
\item $(u,v) \neq (u',v')$ as pairs of elements of $\Gam$.
\end{enumerate}
Given the first condition, the second holds as soon as $u \neq u'$ in $\Gam$,
or equivalently $v \neq v'$. The first condition is an instance of the word
problem and the second is a disequality in $\Gam$. Both have classical
certificates. An equality $w = w'$ in a presented group is certified by writing
$w w'^{-1}$ as a product of conjugates of relators and their inverses in the
free group. A
disequality is certified by a homomorphism to a finite group under which the
two elements have different images.

These two conditions do not require the listed words to represent distinct
elements: if two words of a list coincide in $\Gam$, the sets are smaller than
the lists but the argument is unaffected. Pairwise distinctness is needed only
for the additional statement that $A$ and $B$ have exactly as many elements as
the lists have words. We certify it as well.

\section{The witness and its certificates}\label{sec:certs}

The witness is a pair of lists of reduced words in $a^{\pm1}, b^{\pm1}$: a list
of $32$ words for $A$ and a list of $28$ words for $B$, each word of length at
most six. Both lists are printed in \cref{app:witness}. There are
$32 \cdot 28 = 896$ products.

\subsection{Equality certificates}\label{sec:eqcerts}

For each pair $(u,v)$ of listed words the certificate names a \emph{partner}
$(u',v')$, a different pair of listed words, together with a list of triples
$(w_t, i_t, e_t)$ asserting the identity
\begin{equation}\label{eq:conjrel}
  (uv)(u'v')^{-1}
  \;=\;
  \prod_{t=1}^{m} w_t\, r_{i_t}^{e_t}\, w_t^{-1}
  \qquad\text{in the free group } \Free,
\end{equation}
where $r_1, r_2$ are the relators of \eqref{eq:presentation}, $e_t = \pm 1$,
and the $w_t$ are words. The right-hand side lies in the normal closure of the
relators, so the identity proves $uv = u'v'$ in $\Gam$. Checking it requires
only free reduction of both sides and a comparison of words. Each certificate
is self-contained: it does not refer to any other certified equality.

\begin{example}\label{ex:cert}
Take $u = aba^2$ and $v = b^{-1}$, with partner $u' = \varepsilon$, the empty
word, and $v' = b^2a^{-2}$. Then
$(uv)(u'v')^{-1} = aba^2b^{-1} \cdot a^2b^{-2} = r_1$, so the certificate is
the single triple $(w, i, e) = (\varepsilon, 1, +1)$: the relator $r_1$ itself,
with trivial conjugator. (The generated data number the relators from $0$.)
\end{example}

In $658$ of the $896$ certificates the list of triples is empty: the words
$uv$ and $u'v'$ are already equal in the free group, that is, they have the
same free reduction. The other $238$ certificates use relators, $970$ conjugates in all;
$136$ certificates use a single conjugate, and the largest uses $55$. The
expanded right-hand sides of \eqref{eq:conjrel} have total length $16{,}972$
letters. The search did not try to minimize the use of relators: $710$ of the
$896$ products have some partner that is equal in the free group, so relators
are unavoidable for $186$ products.

\subsection{Disequality certificates}\label{sec:distcerts}

\paragraph{What the theorem uses.}
For each of the $896$ certificates, Lean checks that $u \neq u'$ or
$v \neq v'$ in $\Gam$. The homomorphism $\psi \colon \Gam \to \Zmod{42}$ of
\cref{sec:gamma} is enough for this: in every one of the $896$ certificates,
$\psi(v) \neq \psi(v')$. Lean proves that the images of $a$ and $b$ send both relators to the
identity, so that $\psi$ is well defined on $\Gam$; that the two words have
different images is part of the kernel computation of \cref{sec:lean}, which
accepts a separation by $\psi$ or by $\varphi$ on either coordinate.

\paragraph{Pairwise distinctness.}
To show in addition that $|A| = 32$ and $|B| = 28$ we certify all
$\binom{32}{2} + \binom{28}{2} = 496 + 378 = 874$ disequalities between
listed words. The homomorphism $\psi$ separates $868$ of these pairs,
including all $378$ pairs from the second list. The six pairs it does not
separate are
\[
\begin{array}{lll}
 \{b^2,\; ab^2a^{-1}\}, & \{a^{-1}ba^{-1},\; b^2a^{-2}b^{-1}\}, & \{a^{-1}b^2,\; ab^2a^{-2}\}, \\[2pt]
 \{aba^2,\; a^{-1}ba^{-1}b^2\}, & \{a^{-1}ba^{-1}b,\; b^2a^{-2}\}, & \{a^{-1}b^3,\; a^2ba^2\}.
\end{array}
\]
A second homomorphism $\varphi$ from $\Gam$ to the symmetric group on
$\{0,1,2,3\}$, given by
\[
  \varphi(a) = (0\;1\;2), \qquad \varphi(b) = (1\;3\;2),
\]
sends both relators to the identity and separates each of these six pairs.
The images of $a$ and $b$ are $3$-cycles with different fixed points, so the
image of $\varphi$ is the alternating group $\Afour$; this is consistent with
\eqref{eq:ab}, since the abelianization $\Zmod{3}$ of $\Afour$ is a quotient of
$\Zmod{42}$. All of these verifications are computations in $\Zmod{42}$ and in
the symmetric group on four points, and Lean performs them by \code{decide}.

\subsection{Remarks on the design}

Every true equality in a presented group has a certificate of the form
\eqref{eq:conjrel}, so for equalities the only difficulty is to find one; this
is the purpose of the proof-producing enumeration of \cref{sec:engine}.
Separation by finite quotients is not a complete method for finitely presented
groups in general, which need not be residually finite. For $\Gam$ it is
complete in principle: as an arithmetic lattice, $\Gam$ is a finitely
generated group that is linear over a field (necessarily of positive
characteristic, since $\Gam$ is not linear over $\mathbb{C}$), and finitely
generated linear groups over any field are residually finite by Malcev's
theorem. Like the rest of the description of $\Gam$, this rests
on~\cite{gardamslides2} and is not formalized. What is particular to this
witness is that two very small quotients separate all of its words, and that
$\psi$ alone separates every product from its partner.

\section{The Lean development}\label{sec:lean}

The formalization uses Lean \code{v4.30.0} and Mathlib at commit
\code{c5ea003}. It has four proof modules, $593$ lines in all, two short audit
files, and generated data. \code{Core} defines $\Gam$ as a Mathlib
\code{PresentedGroup} on the two relators of \eqref{eq:presentation}, the
four-letter alphabet, the evaluation map from words to $\Gam$, and the
certificate data types. \code{Reduce} implements free reduction and proves the
soundness of the equality checker: if the expansion of the triples freely
reduces to the same word as $(uv)(u'v')^{-1}$, then $uv = u'v'$ in $\Gam$. It
also proves that the relators written as letter words agree with the relators
of the presentation, so no informal identification is used. \code{Distinct}
constructs $\psi$ and $\varphi$ through \code{PresentedGroup.toGroup}, proves
that a difference of images implies a disequality in $\Gam$, and proves
pairwise distinctness of the two lists. \code{NonUP} checks, in one kernel
computation over the $896$ certificates, that each partner is a pair of listed
words, that the equality certificate is valid, and that $\psi$ or $\varphi$
separates $u$ from $u'$ or $v$ from $v'$; it also checks that the certificates
cover every pair of listed words. From this it derives the result:

\begin{lstlisting}[language=lean4,mathescape=true]
theorem gamma_not_uniqueProds : $\lnot$ UniqueProds Gamma
\end{lstlisting}

\noindent Here \code{Gamma} is the Lean name of $\Gam$.

Pairwise distinctness enters the formal proof only through the cardinalities
$|A| = 32$ and $|B| = 28$, which are used to show that the two sets are
nonempty. The proof as built therefore depends on $\varphi$ as well as on
$\psi$; a direct proof of nonemptiness would remove this dependence, and
mathematically $\psi$ alone suffices.

\paragraph{Generated data.}
The witness words and the $896$ partner certificates are in two files of the
directory \code{Generated/}, written by a code generator from the output of
the search program and never edited by hand. The two homomorphisms are written
directly in \code{Distinct}. Lean checks every certificate from the generated
files by kernel reduction, so an error in the generator or in the search
program can cause the build to fail but cannot cause a false theorem to be
accepted. This is the usual division of labor in verified SAT
checking~\cite{heule2017lrat,lammich2020lrat}.

\paragraph{The axiom gate.}
A kernel-checked theorem may still depend on unexpected axioms. The script
\code{scripts/axiom\_gate.sh}, run after the build and in continuous
integration, audits all $122$ theorem constants in the compiled library,
including those that Lean generates automatically, and fails on any occurrence
of \code{sorry} or \code{native\_decide} and on any axiom outside
$\axiomClosure$. For the main theorem it prints (wrapped here):

\begin{lstlisting}[language={}]
'Karanos.gamma_not_uniqueProds' depends on axioms:
  [propext, Classical.choice, Quot.sound]
\end{lstlisting}

No compiled evaluation is trusted: every \code{decide} in the development is a
kernel reduction.

\section{The search program}\label{sec:engine}

The program that found the witness is not trusted, and nothing in the proof
depends on how it works.

\paragraph{Proof-producing coset enumeration.}
The program builds partial quotients of balls of reduced words by
Todd--Coxeter enumeration~\cite{toddcoxeter1936}, modified so that every
identification of two cosets is recorded together with the relator instance
that forced it. The closure computed on a ball therefore gives, for any two words of the ball
that it identifies, a derivation of their equality.

\paragraph{Extraction.}
A second component converts each required derivation into the flat form
\eqref{eq:conjrel}. Each identification forced directly by a relator contributes
one conjugated relator; identifications that follow from earlier ones
contribute none. Before transcription to Lean, the program checks every certificate by
the same free-reduction computation that Lean later performs.

\paragraph{Search.}
The witness was found by a search over subsets of the enumerated balls, using
a SAT encoding of the condition that no product is uniquely represented.
Homomorphisms to small finite groups were found by a separate search. The
search log is in the repository.

\section{Results}\label{sec:results}

\begin{table}[ht]
  \centering
  \caption{Summary of the development.}
  \label{tab:summary}
  \small
  \begin{tabular}{@{}lr@{}}
    \toprule
    witness sizes $|A|$, $|B|$ & $32$, $28$ \\
    maximum word length in the witness & $6$ \\
    product coincidences certified & $896$ \\
    \quad certified by an identity in the free group & $658$ \\
    \quad certified with relators & $238$ \\
    conjugated relators used (total; largest certificate) & $970$; $55$ \\
    certificates in which $\psi$ separates $v$ from $v'$ & $896$ \\
    pairwise disequalities certified & $874$ \\
    \quad separated by $\psi \colon \Gam \to \Zmod{42}$ & $868$ \\
    \quad separated only by $\varphi \colon \Gam \to \Afour$ & $6$ \\
    theorem constants audited by the axiom gate & $122$ \\
    axioms permitted & $\axiomClosure$ \\
    occurrences of \code{sorry} / \code{native\_decide} & $0$ / $0$ \\
    \bottomrule
  \end{tabular}
\end{table}

\Cref{tab:summary} summarizes the development. As a check that is independent
of Lean and of the search program, the repository contains a second certificate
checker, a Python script that uses only the standard library
(\code{recheck\_certificates.py}, in the directory \code{scripts}). It reads the two generated data files; the relators and the
two homomorphisms are restated in the script. It confirms that the $896$
certificates cover every pair of listed words, that no pair is its own
partner, that every identity \eqref{eq:conjrel} holds under free reduction,
and that the two homomorphisms send the relators to the identity and separate
all $874$ pairs, and that $\psi$ alone separates $v$ from $v'$ in every
certificate. It reproduces the figures of \cref{tab:summary} that concern the
certificates. Lean itself checks only the weaker disjunction described in
\cref{sec:lean}.

\section{Limitations}\label{sec:limits}

\begin{enumerate}[leftmargin=2em]
\item The formal theorem concerns the group presented by
  \eqref{eq:presentation}. That this group is torsion-free, that it is an
  $\Atwo$ lattice, and that it has property~(T) are not formalized here. The
  relevance of the theorem to Kaplansky's conjecture depends on
  torsion-freeness, which we take from Gardam.
\item Lean checks that $\psi$ and $\varphi$ are well defined and that they
  separate the required pairs. It does not check that $\Gab$ is isomorphic to
  $\Zmod{42}$ or that the image of $\varphi$ is $\Afour$; these are the short
  computations given in \cref{sec:gamma,sec:distcerts}, and the proof does not
  use them.
\item We make no claim that the witness is minimal in size or in word length.
\item The search program is tested but not verified. This does not affect the
  theorem, which depends only on the generated data that Lean re-checks.
\end{enumerate}

\section{Related work}\label{sec:related}

\paragraph{Groups without unique products.}
The known examples are listed in \cref{sec:background}. Recent computational
work studies sets without unique products quantitatively. Dietrich, Lee, Nies
and Vinyals~\cite{dlnv2026} report computational experiments on the
trivial-units and unique-product properties, including a new candidate group.
Tabei~\cite{tabei2026a,tabei2026b,tabei2026c} studies sets without unique
products inside balls of the Promislow group and in the groups of Nielsen and
Soelberg~\cite{nielsensoelberg2024}, and the support of units in $\Ftwo[P]$;
negative results there are certified by constraint and SAT solvers, in part
with DRAT proofs, and positive witnesses by independent scripts. None of this
work uses a proof assistant, and none of it concerns the group $\Gam$.

\paragraph{Formalized results on group rings.}
Gadgil and Tadipatri~\cite{gadgiltadipatri2024} formalized Gardam's
counterexample to the unit conjecture in Lean~4. A group with unique products
has only trivial units in its group rings over fields, so their theorem implies, by
a standard argument that is not formalized there as far as we can tell, that
$P$ does not have unique products. Because we work from a
presentation of the group, our development is organized around certificates
for the word problem.

\paragraph{Earlier work with the same axiom gate.}
We used the same organization and the same axiom gate in two earlier papers:
exclusion results for the Erd\H{o}s--Selfridge covering
problem~\cite{miansiddique2026}, and certificates for the geometric part of the
lower bound in the minimum Kochen--Specker problem~\cite{siddiquemian2026}.

\section{Conclusion}\label{sec:conclusion}

Gardam's announced theorem that $\Gam$ does not have unique products now has a
proof checked by the Lean kernel against Mathlib's definition, with an explicit
witness and certificates that can be checked independently of Lean.

\paragraph{Artifact availability.}
The Lean development, the search program, the certificate data and the
independent checker are available under the Apache~2.0 license at
\begin{center}\url{https://github.com/ibrahimmian36/Karanos}\end{center}
The version described here is commit \code{fa39e9a}.

\needspace{10\baselineskip}
\section{Use of AI tools}\label{sec:ai}

A large language model (LLM) was used extensively in this work: Claude
(Anthropic), between July and September 2026. Under the authors' direction the
LLM wrote most of the Lean formalization and of the Python search program, and
it generated much of the text of this paper, which the authors reviewed and
revised. The authors take full intellectual responsibility for the entire
content of the paper, including the text, the references and every claim. The
main theorem is checked by the Lean kernel and does not depend on the
correctness of LLM-written code.

\paragraph{Acknowledgments.}
We thank the maintainers of Lean and Mathlib.

\appendix

\needspace{20\baselineskip}
\section{The witness}\label{app:witness}

The two lists below are the witness, copied mechanically from the generated
file \code{Witness.lean} and printed in the order of that file. Words are written in the letters
\code{a}, \code{b}, \code{A}, \code{B}, where \code{A} $= a^{-1}$ and
\code{B} $= b^{-1}$, and $\varepsilon$ is the empty word. All words are freely
reduced.

\medskip
\noindent\textbf{The set $A$ ($32$ words).}
\begin{center}\small
\begin{tabular}{@{}llllllll@{}}
$\varepsilon$ & \texttt{a} & \texttt{A} & \texttt{b} & \texttt{ab} & \texttt{aB} & \texttt{Ab} & \texttt{bb} \\
\texttt{aba} & \texttt{abb} & \texttt{aBB} & \texttt{AbA} & \texttt{Abb} & \texttt{abaa} & \texttt{abab} & \texttt{abbA} \\
\texttt{abbb} & \texttt{AbAb} & \texttt{Abbb} & \texttt{aabaa} & \texttt{abaab} & \texttt{bbAA} & \texttt{ababb} & \texttt{abbAA} \\
\texttt{abbAb} & \texttt{abbbb} & \texttt{AbAA} & \texttt{AbAbb} & \texttt{bbAAB} & \texttt{abaabb} & \texttt{abbAbb} & \texttt{AbAAb} \\
\end{tabular}
\end{center}

\noindent\textbf{The set $B$ ($28$ words).}
\begin{center}\small
\begin{tabular}{@{}llllllll@{}}
$\varepsilon$ & \texttt{A} & \texttt{b} & \texttt{B} & \texttt{ab} & \texttt{AA} & \texttt{bA} & \texttt{bb} \\
\texttt{BA} & \texttt{bAA} & \texttt{bbA} & \texttt{bbb} & \texttt{Bab} & \texttt{BBAB} & \texttt{bbAA} & \texttt{abbab} \\
\texttt{AbAA} & \texttt{BAbAA} & \texttt{BBAAA} & \texttt{BBABA} & \texttt{BBBAB} & \texttt{abbbbb} & \texttt{babbab} & \texttt{bbAbAA} \\
\texttt{BBBAAA} & \texttt{BBBABA} & \texttt{BBABab} & \texttt{bAbAA} &  &  &  &  \\
\end{tabular}
\end{center}


\begin{thebibliography}{99}

\bibitem{agol2013}
I.~Agol.
\newblock The virtual {H}aken conjecture.
\newblock With an appendix by I.~Agol, D.~Groves and J.~Manning.
\newblock \emph{Documenta Mathematica}, 18:1045--1087, 2013.

\bibitem{arzhantsevasteenbock2023}
G.~Arzhantseva and M.~Steenbock.
\newblock Rips construction without unique product.
\newblock \emph{Pacific Journal of Mathematics}, 322(1):1--9, 2023.

\bibitem{bhv2008}
B.~Bekka, P.~de~la Harpe and A.~Valette.
\newblock \emph{Kazhdan's Property (T)}.
\newblock Cambridge University Press, 2008.

\bibitem{carter2014}
W.~Carter.
\newblock New examples of torsion-free non-unique product groups.
\newblock \emph{Journal of Group Theory}, 17(3):445--464, 2014.

\bibitem{cmsz1993}
D.~I. Cartwright, A.~M. Mantero, T.~Steger and A.~Zappa.
\newblock Groups acting simply transitively on the vertices of a building of
  type $\tilde{A}_2$, {I}.
\newblock \emph{Geometriae Dedicata}, 47(2):143--166, 1993.

\bibitem{cms1994}
D.~I. Cartwright, W.~M{\l}otkowski and T.~Steger.
\newblock Property ({T}) and $\tilde{A}_2$ groups.
\newblock \emph{Annales de l'Institut Fourier}, 44(1):213--248, 1994.

\bibitem{heule2017lrat}
L.~Cruz-Filipe, M.~J.~H. Heule, W.~A. Hunt~Jr., M.~Kaufmann and
  P.~Schneider-Kamp.
\newblock Efficient certified {RAT} verification.
\newblock In \emph{Automated Deduction -- CADE 26}, LNCS 10395, pages 220--236.
  Springer, 2017.

\bibitem{dlnv2026}
H.~Dietrich, M.~Lee, A.~Nies and M.~Vinyals.
\newblock On the trivial units property and the unique product property.
\newblock arXiv:2603.22640, 2026.

\bibitem{fishersanchezperalta}
S.~P. Fisher and P.~S\'anchez-Peralta.
\newblock Division rings for group algebras of virtually compact special groups
  and 3-manifold groups.
\newblock \emph{Journal of Combinatorial Algebra}, 10(1/2):153--193, 2026.
  arXiv:2303.08165.

\bibitem{gadgiltadipatri2024}
S.~Gadgil and A.~R. Tadipatri.
\newblock Formalizing {G}iles {G}ardam's disproof of {K}aplansky's unit
  conjecture.
\newblock In \emph{Proceedings of the 13th ACM SIGPLAN International Conference
  on Certified Programs and Proofs (CPP 2024)}, pages 177--189. ACM, 2024.

\bibitem{gardam2021}
G.~Gardam.
\newblock A counterexample to the unit conjecture for group rings.
\newblock \emph{Annals of Mathematics}, 194(3):967--979, 2021.

\bibitem{gardamslides1}
G.~Gardam.
\newblock Kaplansky's conjectures.
\newblock Slides of a lecture at the Global Noncommutative Geometry Seminar,
  17 September 2021.
\newblock \url{https://www.gilesgardam.com/slides/gncg.pdf}.

\bibitem{gardamslides2}
G.~Gardam.
\newblock The {K}aplansky conjectures.
\newblock Slides of a talk in the Geometry and Topology Online seminar,
  University of Warwick, 25 November 2021.
\newblock \url{https://sschleimer.warwick.ac.uk/Seminar/Talks/2021-11-25gardam.pdf}.

\bibitem{gms2015}
D.~Gruber, A.~Martin and M.~Steenbock.
\newblock Finite index subgroups without unique product in graphical small
  cancellation groups.
\newblock \emph{Bulletin of the London Mathematical Society}, 47(4):631--638,
  2015.

\bibitem{kaplansky1970}
I.~Kaplansky.
\newblock ``{P}roblems in the theory of rings'' revisited.
\newblock \emph{American Mathematical Monthly}, 77(5):445--454, 1970.

\bibitem{lammich2020lrat}
P.~Lammich.
\newblock Efficient verified ({UN}){SAT} certificate checking.
\newblock \emph{Journal of Automated Reasoning}, 64(3):513--532, 2020.

\bibitem{mathlib2020}
{The mathlib Community}.
\newblock The {L}ean mathematical library.
\newblock In \emph{Proceedings of the 9th ACM SIGPLAN International Conference
  on Certified Programs and Proofs (CPP 2020)}, pages 367--381. ACM, 2020.

\bibitem{miansiddique2026}
I.~Mian and S.~Siddique.
\newblock Kernel-checked exclusions for the {E}rd\H{o}s--{S}elfridge odd
  covering problem: any odd covering of $\mathbb{Z}$ has lcm exceeding 10000.
\newblock arXiv:2607.25628, 2026.

\bibitem{moura2021lean}
L.~de~Moura and S.~Ullrich.
\newblock The {L}ean 4 theorem prover and programming language.
\newblock In \emph{Automated Deduction -- CADE 28}, LNCS 12699, pages 625--635.
  Springer, 2021.

\bibitem{nibloreeves1997}
G.~Niblo and L.~Reeves.
\newblock Groups acting on {CAT}(0) cube complexes.
\newblock \emph{Geometry \& Topology}, 1:1--7, 1997.

\bibitem{nielsensoelberg2024}
P.~P. Nielsen and L.~Soelberg.
\newblock Small sets without unique products in torsion-free groups.
\newblock \emph{Journal of Algebra and Its Applications}, 23(8), article
  2550050, 2024.

\bibitem{passman1977}
D.~S. Passman.
\newblock \emph{The Algebraic Structure of Group Rings}.
\newblock Wiley-Interscience, 1977.

\bibitem{promislow1988}
S.~D. Promislow.
\newblock A simple example of a torsion-free, non unique product group.
\newblock \emph{Bulletin of the London Mathematical Society}, 20(4):302--304,
  1988.

\bibitem{ripssegev1987}
E.~Rips and Y.~Segev.
\newblock Torsion-free group without unique product property.
\newblock \emph{Journal of Algebra}, 108(1):116--126, 1987.

\bibitem{siddiquemian2026}
S.~Siddique and I.~Mian.
\newblock Machine-checked certificates for the geometric half of the minimum
  {K}ochen--{S}pecker bound.
\newblock arXiv:2607.26413, 2026.

\bibitem{soelberg2018}
L.~J. Soelberg.
\newblock Finding torsion-free groups which do not have the unique product
  property.
\newblock Master's thesis, Brigham Young University, 2018.
\newblock \url{https://scholarsarchive.byu.edu/etd/6932/}.

\bibitem{steenbock2015}
M.~Steenbock.
\newblock Rips--{S}egev torsion-free groups without the unique product
  property.
\newblock \emph{Journal of Algebra}, 438:337--378, 2015.

\bibitem{tabei2026a}
M.~Tabei.
\newblock Least sizes of non-unique-product sets: the {P}romislow group and a
  {H}eisenberg-type candidate.
\newblock arXiv:2607.18346, 2026.

\bibitem{tabei2026b}
M.~Tabei.
\newblock The quantitative non-unique-product landscape at the global minimum:
  the {N}ielsen--{S}oelberg groups.
\newblock arXiv:2607.19687, 2026.

\bibitem{tabei2026c}
M.~Tabei.
\newblock Localizing the {G}ardam unit: the support geometry of units in
  $\mathbb{F}_2[P]$ and its non-unique-product relatives.
\newblock arXiv:2609.17559, 2026.

\bibitem{toddcoxeter1936}
J.~A. Todd and H.~S.~M. Coxeter.
\newblock A practical method for enumerating cosets of a finite abstract group.
\newblock \emph{Proceedings of the Edinburgh Mathematical Society}, 5(1):26--34,
  1936.

\bibitem{wise2021}
D.~T. Wise.
\newblock \emph{The Structure of Groups with a Quasiconvex Hierarchy}.
\newblock Annals of Mathematics Studies 209. Princeton University Press, 2021.

\end{thebibliography}
\end{document}